\documentclass[a4paper,12pt]{article}

\usepackage{amsmath}
\usepackage{amssymb}
\usepackage{epsfig}
\usepackage{subfigure}
\usepackage{float}
\usepackage{url}
\usepackage{enumerate}
\usepackage[T1]{fontenc}

\newtheorem{theo}{thm}[section]
\newtheorem{conjecture}[theo]{Conjecture}
\newtheorem{corollary}[theo]{Corollary}
\newtheorem{definition}[theo]{Definition}
\newtheorem{lemma}[theo]{Lemma}
\newtheorem{notation}[theo]{Notation}
\newtheorem{proposition}[theo]{Proposition}
\newtheorem{question}[theo]{Question}
\newtheorem{remark}[theo]{Remark}
\newtheorem{theorem}[theo]{Theorem}

\def\C{\mathbb{C}}
\def\N{\mathbb{N}}
\def\Q{\mathbb{Q}}
\def\R{\mathbb{R}}
\def\Z{\mathbb{Z}}

\def\P{\mathcal{P}}
\def\PP{\text{\rm \bf P}}

\def\card{\textrm{\rm card}}
\def\ed{\textrm{\rm Ed}}
\def\l{\ell}
\def\endproof{\hfill  $\square$ \bigskip}

\def\pmm{p_{i^-}^-}
\def\qmm{q_{i^-}^-}
\def\pp{p_{i^-}^-}
\def\qq{q_{i^-}^-}
\def\pgcd{\textrm{\rm gcd}}
\def\proof{\noindent {\bf Proof. }}
\def\splice{Eisenbud and Neumann }
\def\sp{\textrm{\rm Sp}}
\def\spp{\textrm{\rm Spp}}
\def\sppa{\textrm{\rm Sppa}}

\title{The Hertling conjecture in dimension 2}
\author{Thomas Br\'elivet}
\date{}

\begin{document}

\maketitle
\def\thefootnote{\relax}
\footnote{
 {\bf \hskip -4ex AMS classification (2000): }14B05, 14H20, 52C05, 
        }
Dedicated to the memory of my father who died during the preparation of this manuscript.
\begin{abstract}
We consider an isolated plane curve singularity and its associated \splice diagram. We give an algorithm to compute the maximal spectral value on the diagram and we show that the singularity is topologically equivalent to another singularity for which the maximal spectral value is given by the point $(1,1)$ in the plane of the Newton polygon. From the almost additivity on the splice components of the diagram we compute the sum of the square of the spectral values. This formula with the previous result on the maximal spectral value give us the Hertling conjecture as a corollary.
\end{abstract}

\def\abstractname{Résumé}

\begin{abstract}
On considère une singularité isolée de courbe ainsi que son diagramme de Eisenbud et Neumann. On donne un algorithme pour calculer la valeur spectrale maximale sur le diagramme et on montre que la singularité est topologiquement équivalente à une autre telle que sa valeur spectrale maximale est donnée par le point $(1,1)$ dans le plan du polygone de Newton. De la presque additivité sur les composantes de splice du diagramme, on calcule la somme des carrés des valeurs spectrales. Cette formule ainsi que les résultats précédents sur la valeur spectrale maximale nous donne la conjecture de Hertling comme corollaire.
\end{abstract}

\tableofcontents

\section{Introduction}

Some 25 years ago, Steenbrink has defined the spectrum of an isolated hypersurface singularity in \cite{St1} and then Steenbrink himself, Varchenko \cite{V} and others have obtained very interesting results motivated mainly by one conjecture made by Arnold, see \cite{AGV} and \cite{St2} for more details.

The spectrum is a collection of rational numbers between $-1$ and $n$, where $n+1$ denotes the dimension of the ambient space and it is symmetric around $(n-1)/2$.

The variance measures the distribution of these numbers with respect to the middle point and is defined by 

$$V=\frac{1}{\mu}\sum_{i=1}^{\mu}\left( \alpha_i-\frac{n-1}{2} \right)^2$$
where $\alpha_1+ \cdots +\alpha_\mu$ as an element of $\N^{(\Q)}$ is the spectrum with $\alpha_1 \leq \cdots \leq \alpha_\mu$. 

It came as a great surprise when Hertling, at  the Summer Institute on Singularities, Newton Institute, Cambridge 2000, proposed the following conjecture.

\begin{conjecture}\label{conjecture : Hertling}
For any isolated hypersurface singularity
$$V \leq \frac{\alpha_\mu-\alpha_1}{12}.$$
\end{conjecture}

This conjecture was supported at the time by the case of weighted homogeneous singularities where one has in fact an equality (see \cite{H} for a conceptual proof involving Frobenius manifolds and \cite{Di} for a high school proof based on some formulas in \cite{St1}) as well as by inspection through Arnold's lists of singularities.

Soon after this, M. Saito (see \cite{S2}) showed that Conjecture \ref{conjecture : Hertling} holds for all irreducible plane curves singularities. 
In 2002 it has been proved by the author that Conjecture \ref{conjecture : Hertling} also holds for all non-degenerated and commode plane curves singularities.

Here we prove the Conjecture \ref{conjecture : Hertling} for all isolated plane curve singularities, see Corollary \ref{corollary : Hertling conjecture}.

In section \ref{section : Definition of the spectral pairs of a fibered splice diagram} we recall how to compute the spectral pairs of an isolated curve singularity.

Note that in every dimension the spectrum of a Newton non-degenerated singularity is known from the Newton polyhedron by Steenbrink \cite{St1}, Khovanski\u\i\ and Varchenko \cite{Khovanskii Varchenko}.

In section \ref{setcion : Spectral pairs of a parallelogram} we give a geometric description in terms of Newton polygons, generalizing the well known situation of the Newton Non-degenerated case.

In section \ref{section: the non-degenerated and commode case} we recall the formulas in the non-degenrated case.

In section \ref{section : Additivity of the spectral pairs} we construct an application from the set of isolated singularities plane curve in the free group generated by the commode Newton polygons.

In section \ref{section : The maximal spectral value} we give an algorithm to compute the maximal spectral value and we prove the Theorem \ref{Theorem : The spectral maximal value}.

In every dimension we already know that the multiplicity of the minimal spectral value is one, see \cite{S1}. We Thank Antoine Douai for this reference.

Finally the section \ref{section : A formula for the variance of the spectrum} is the core of the proof of the conjecture and gives us an expression of $6S-\mu\alpha_\mu$ ($S$ is the sum of the square of the spectral values and $\alpha_\mu$ the maximal spectral value) as a linear combination with strictly negative coefficients of the determinants of the \splice diagram representing the link of $f$. The formula gives the Hertling conjecture as a corollary.

\section{Definition of the spectral pairs of a fibered \splice diagram}\label{section : Definition of the spectral pairs of a fibered splice diagram}

In this section we define a notion of spectral pairs associated to a fibered multilink, see \cite{SSS} and \cite{C}. Let $L$ be a fibered multilink and $(\Gamma,*)$ a rooted \splice diagram with non zero determinant representing $L$ (for instance minimal). 
See \cite{Eisenbud Neumann}, \cite{CP} for a complete introduction of \splice diagram and \cite{Neumann} for a rapid and historic introduction.

Let $V$ the set of vertices, $\ed$ the set of edges, $A$ the set of arrows and $R$ the set of rupture vertices (vertices such that the number of incident edges is greater or equal than $2$) of $\Gamma$.

\begin{figure}[H]
\begin{center}
\input{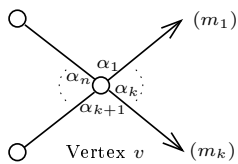}
\caption{Neighborhood of $v$}\label{Figure : Neighborhood of v}
\end{center}
\end{figure}

Let $v$ be a vertex of $\Gamma$.
Cut edges joining $v$ and rupture vertices. Replace the edges by arrows with multiplicities such that the multiplicity $m_v$ of $v$ doesn't change (see Figure \ref{Figure : Neighborhood of v}).
We define $m_i=0$ for $i=k+1, \cdots, n_v$, $\beta_j$ for $j=1,\cdots, n_v$ such that 
$$
\beta_j \alpha_1 \cdots \hat{\alpha_j} \cdots \alpha_{n_v} \equiv 1 \mod\alpha_j
$$

and 
$$
s_{v,w}=s_{v,j}=\frac{m_j-\beta_jm_v}{\alpha_j}, \textrm{ where } w \textrm{ is the end of the } j^{\textrm{th}} \textrm{edge}.
$$

\begin{enumerate}[(1)]
\item If $v$ is not the root then let $p(v)$ be the predecessor of $v$ given by the first vertex met in the path from $v$ to $*$ (it is unic because $\Gamma$ is a tree),
\item $m_v$ the multiplicity of the vertex $v$, it is the sum over the arrows of the product of all edge weights adjacent to the path from $v$ to the arrows, 
\item $d_v=\gcd(m_v,s_{v,p(v)})$, 
\item $r_v=\gcd(m_v,s_{v,p(v)},j=1,\cdots,n_v)$. 
\end{enumerate}

Then we define some elements of $\Z^{(\Q\times\Z)}$ by 
\begin{enumerate}
\item[]
$\displaystyle
a_v=\sum_{\underset{m_v \nmid s r_v}{0<s<m_v}} \left( -1+\sum_{i=1}^{n_v} \left\lbrace \frac{ss_{v,i}}{m_v} \right\rbrace \right) \left[ \left( \frac{s}{m_v}-1,1 \right) + \left( 1-\frac{s}{m_v},1 \right) \right] , \text{  } v \in V \cup \{ * \},
$
\item[] 
$\displaystyle
b_v=\sum_{0<s<r_v}\left[ \left( -\frac{s}{r_v},2 \right) + \left( \frac{s}{r_v},0 \right) \right] , \text{  } v \in V \cup \{ * \},
$
\item[] 
$\displaystyle
c_v=\sum_{0<s<d_v}\left[ \left( -\frac{s}{d_v},2 \right) + \left( \frac{s}{d_v},0 \right) \right], \text{  } v \in V,
$
\item[] 
$\displaystyle
c'_v=\sum_{0<s<d_v}\left( -\frac{s}{d_v},2 \right) , \text{  } v \in A, 
$
\end{enumerate} 
where $\{ x \}$ for $x \in \R$ means the fractional part of $x$.

\begin{definition}
The spectral pairs of $L$ are defined by
$$
\spp(L)=\sum_{v\in R}a_v+\sum_{v\in R \setminus \{ * \} }(c_v-b_v)-b_*+\sum_{v \in A}c'_v+(|A|-1)(0,1).
$$
and the spectrum of $L$ is defined by the projection on the first factor of the spectral pairs and is denoted by $\sp(L)$.
\end{definition}

\begin{remark}\rm\label{remark : invariance of the spectrum with zero determinants}
The definition of the spectral pairs of a fibered multilink is an invariant of the topology of the complementary of the link and is independant of the choice of the root. It is also independant of the diagram if we not permit zero determinants. The spectrum is independant of the choice of the root and the diagram (even if we accept zero determinants). 
\end{remark}

\begin{theorem}[\cite{SSS}]
Let $f$ be an isolated plane curve singularity and $L_f$ the link associated to $f$. Then $\spp(L_f)=\spp(f)$. 
\end{theorem}

The following proposition shows that $\spp$ is almost additive. The section \ref{section : Additivity of the spectral pairs} will explain that through a factorisation $\spp$ is additive.

\begin{proposition}[\cite{SSS}]
Suppose that the fibered multilink $L$ is the result of splicing the fibered multilinks $L_1$ and $L_2$ along components of multilink multiplicities $m_1$ and $m_2$. Let $d=\gcd(m_1,m_2)$ then
$$
\spp(L)=\spp(L_1)+\spp(L_2)-(0,1)+\sum_{s=1}^{d-1}\left[ \left( \frac{s}{d}, 0 \right) - \left( - \frac{s}{d}, 2 \right) \right]. 
$$
\end{proposition}

We denote by $\alpha_1, \ldots, \alpha_\mu$ the spectral values, where $\mu$ is the number of spectral values counted with multiplicity and $\alpha_1 \leq \cdots \leq \alpha_\mu$.

\section{Geometric description of the spectrum}\label{setcion : Spectral pairs of a parallelogram}

\begin{figure}[H]
  \centering
  \input{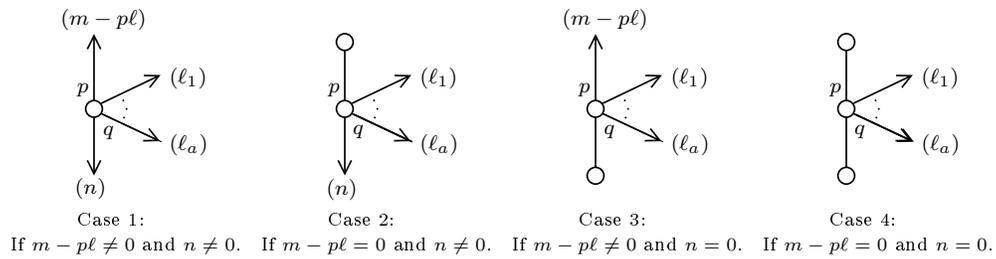}
  \caption{Parallelogram Diagrams}\label{Figure : ElementaryDiagram}
\end{figure}

Let us consider $\Gamma(m,n,p,q,\l_1,\ldots,\l_a)$ the (minimal) diagram defined by the Figure \ref{Figure : ElementaryDiagram} where $m$, $n$ are non negative integers,
$(p,q)$ are coprime positive integers, $\l_1,\ldots,\l_a$ are positive numbers, $\l=\l_1+\cdots+\l_a$, with $m-p\l>0$.

Diagrams of this type are the bricks of the \splice diagrams.

\begin{figure}[H]
  \centering
  \input{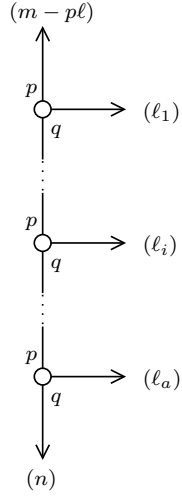}
  \caption{Separated Arrows}\label{Figure : Separated Arrows}
\end{figure}

As we have seen in the Remark \ref{remark : invariance of the spectrum with zero determinants} the spectrum of $\Gamma(m,n,p,q,\l_1,\ldots,\l_a)$ and of the Figure \ref{Figure : Separated Arrows} are equals.

So from now we suppose that $a=1$ and $\l_1=\l$ so the diagram is $\Gamma(m,n,p,q,\l)$ and we compare its spectrum with the spectrum of $\Gamma(p,q,m,n,\underbrace{1,\ldots,1}_{\l \textrm{ times}})$.

Chose $u$ and $v$ such that $up+vq=1$. And note $r=\pgcd(\l,m,n)$, $d=\pgcd(\l,N)$, $d_0=\pgcd(m-p\l,n+q\l)$ and $d_\infty=\pgcd(m,n)$.
Then we have
\begin{multline*}
\spp(\Gamma(m,n,p,q,\l))= \\
\sum_{\substack{0<s<N \\ N \nmid sr}} 
\left( -1 
+ \left\lbrace \frac{s\l}{N} \right\rbrace 
+ \left\lbrace \frac{s(um-vn-\l)}{N} \right\rbrace 
+ \left\lbrace \frac{s(vn-um)}{N} \right\rbrace 
\right)
\left(
\pm \left( \frac{s}{N}-1 \right),1
\right) \\
- \sum_{0<s<r} \left[ \left( - \frac{s}{r},2 \right) + \left( \frac{s}{r},0 \right) \right]
+ \sum_{0<s<d} \left( - \frac{s}{d},2 \right) 
+ B(\Gamma(m,n,p,q,\l))
+ 2 (0,1),
\end{multline*}
and
\begin{multline*}
\spp(\Gamma(p,q,m,n,\underbrace{1,\ldots,1}_{\l \textrm{ times}}))= \\
\sum_{\substack{0<s<N}} 
\left( -1 
+ \frac{s\l}{N} 
+ \left\lbrace \frac{s(um-vn-\l)}{N} \right\rbrace 
+ \left\lbrace \frac{s(vn-um)}{N} \right\rbrace 
\right)
\left(
\pm \left( \frac{s}{N}-1 \right),1
\right) \\
+ B(\Gamma(p,q,m,n,\underbrace{1,\ldots,1}_{\l \textrm{ times}}))
+ (t+1) (0,1),
\end{multline*}
where

$$
B(\Gamma(m,n,p,q,\l))=
\left\lbrace
\begin{array}{lc}
\sum_{s=1}^{d_0-1} \left( - \frac{s}{d_0},2 \right) + \sum_{s=1}^{d_\infty-1} \left( - \frac{s}{d_1},2 \right) & \textrm{Case 1,} \\
\sum_{s=1}^{d_1-1} \left( - \frac{s}{d_\infty},2 \right) & \textrm{Case 2,} \\
\sum_{s=1}^{d_0-1} \left( - \frac{s}{d_0},2 \right) & \textrm{Case 3,} \\
0 & \textrm{Case 4,} 
\end{array}
\right.
$$

and
$B(\Gamma(p,q,m,n,\underbrace{1,\ldots,1}_{\l \textrm{ times}}))=B(\Gamma(m,n,p,q,\l))$.

As we will see the following Lemma is very usefull to understand the geometry of the spectrum.

\begin{lemma}
Let $x$ and $y$ be real numbers such that $x<y$. Then we have
$$
\card(]x,y[ \cap \Z ) = y - x + \left\lbrace x \right\rbrace + \left\lbrace -y \right\rbrace -1. 
$$
\end{lemma}

Let $\phi$ the linear map defined by
$$
\begin{array}{cccc}
\phi : & \R^2  & \longrightarrow & \R \\   
       & (x,y) & \longmapsto     & \frac{qx+py}{qm+pn}.
\end{array}
$$

We then have a generalisation relative to a cone.

\begin{lemma}
Let $s$ in $\Z$, $(m_0,n_0),(m_1,n_1)\in\left(\N^* \right)^2$, $\Q$ linearly independants, $K$ the open cone in $\R_+^2$ generated by $(m_0,n_0)$ and $(m_1,n_1)$. Then we have
\begin{multline*}
\card \left( K \cap \phi^{-1} \left( \frac{s}{p_1n_1+q_1m_1} \right) \cap \N^2 \right) = \\
-1+\frac{k_1}{p_1n_1+q_1m_1}s + \left\lbrace \frac{um_0-vn_0}{p_1n_1+q_1m_1}s \right\rbrace + \left\lbrace \frac{vn_1-um_1}{p_1n_1+q_1m_1}s \right\rbrace.
\end{multline*}
\end{lemma}

Now we are ready to write a geometric interpretation of the spectrum in the non-degenerated case.

\begin{proposition}
\begin{multline*}
\spp(\Gamma(p,q,m,n,\underbrace{1,\ldots,1}_{\l \textrm{ times}}))= \\
\sum_{(m,n) \in P} \left( 1-\phi(m,n),1 \right) + B(\Gamma(p,q,m,n,\underbrace{1,\ldots,1}_{\l \textrm{ times}})) + (t+1) \left( 0,1 \right),
\end{multline*}
where $P$ is the open parallelogram generated by $(m,n)$ and $(m-p\l,n+q\l)$.
\end{proposition}

Now the following formula permits us to understand the geometry in the non degenerated case. 

\begin{multline*}
\spp(\Gamma(m,n,p,q,\l))
= 
\spp(\Gamma(p,q,m,n,\underbrace{1,\ldots,1}_{\l \textrm{ times}})) \\
- \sum_{0<s<N} \left( \frac{s\l}{N} - \left\lbrace \frac{s\l}{N} \right\rbrace \right) \left( \pm \left( \frac{s}{N}-1 \right),1 \right)
+ \sum_{\substack{0<s<N \\ N / sr}} \left( \pm \left( \frac{s}{N}-1 \right),1 \right) \\
- \sum_{0<s<r} \left[ \left( - \frac{s}{r},2 \right) + \left( \frac{s}{r},0 \right) \right]
+ \sum_{0<s<d} \left( - \frac{s}{d},2 \right) 
- (t-1) (0,1).
\end{multline*}

\begin{figure}[H]
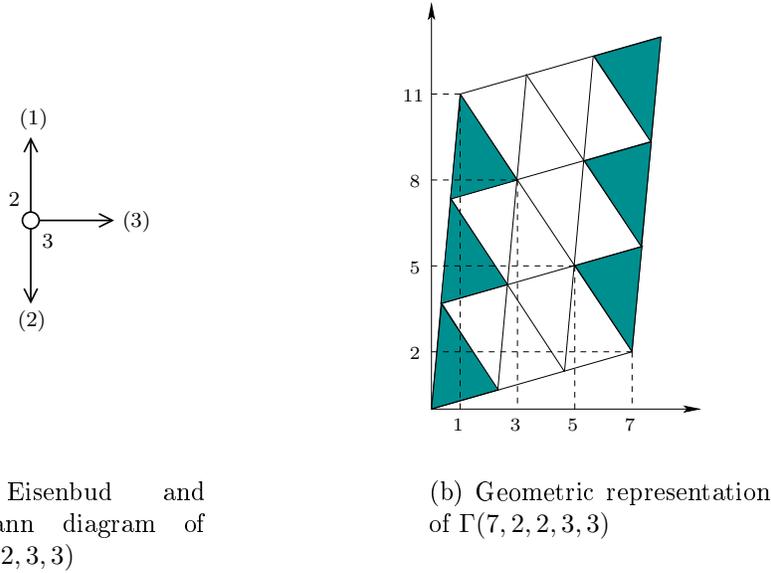

\centering
\begin{center}
  \subfigure[\splice diagram of $\Gamma(7,2,2,3,3)$]{\input{degamma.pstex_t}\label{Figure: splice diagram example 72233}} \hspace{2cm}
  \subfigure[Geometric representation of $\Gamma(7,2,2,3,3)$]{\input{deg.pstex_t}\label{Figure: geometric representation of 72233}} 
\end{center}
\caption{Example of geometric representation of the spectrum}
\label{Figure: Example of geometric representation of the spectrum}
\end{figure}

To avoid long explanation we give an example of geometric representation of $\Gamma(7,2,2,3,3)$ in the Figure \ref{Figure: Example of geometric representation of the spectrum} in order to understand the geometry of the spectrum when it is degenerated with respect to the Newton polygon.

\section{The non-degenerated and commode case}
\label{section: the non-degenerated and commode case}

Let $f : (\C^2,0) \rightarrow (\C,0)$ be an isolated singularity of curve. 
We suppose that $f$ is Newton non-degenerated and commode.

\begin{figure}[H]
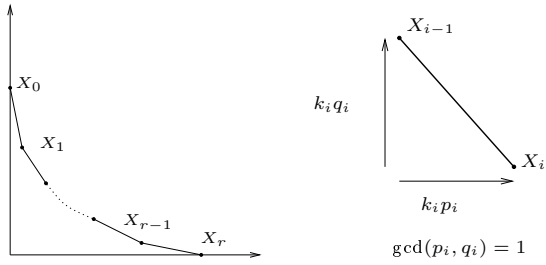

\begin{center}
  \subfigure[Local case]{\input{polygonelocal.pstex_t}\label{Figure: Newton polygons Local}} \quad
  \subfigure[One face]{\input{kipiqi.pstex_t}\label{Figure: Newton polygons kipiqi}} 
\end{center}
\caption{Newton polygon}
\label{Figure: Newton polygons}
\end{figure}

The Newton polygon is given by the Figure \ref{Figure: Newton polygons Local} 

and the \splice diagram is given by the Figure \ref{Figure : diagram in the Newton non-degenerated and commode case}.
\begin{figure}[H]
  \centering
  \input{nondeg.pstex_t}
  \caption{\splice diagram}\label{Figure : diagram in the Newton non-degenerated and commode case}
\end{figure}

The correspondance between $(X_i=(m_i,n_i))_{0 \leq i \leq r}$ ($r$ is the number of faces of the Newton polygon of $f$) and $(p_i,q_i,k_i)_{1 \leq i \leq r}$ is given by
$$
\begin{array}{l}
m_i=k_1p_1+ \cdots + k_ip_i,  \\
n_i=k_{i+1}q_{i+1}+ \cdots + k_rq_r,
\end{array}
$$
for $0 \leq i \leq r$
and
$$
\begin{array}{l}
k_i=\gcd(m_i-m_{i-1},n_{i-1}-n_i), \\
p_i=(m_i-m_{i-1})/k_i, \\
q_i=(n_{i-1}-n_i)/k_i.
\end{array}
$$
for $1 \leq i \leq r$.

It is usefull to introduce two new points $X_{-1}=(1,1)$, $X_0=(1,1)$ and the following notations:
\begin{enumerate}
\item[] $d_i=\pgcd(m_i,n_i)$, $1 \leq i \leq r$,
\item[] $A_{i,j}=m_in_j-m_jn_i$, $-1 \leq i,j \leq r+1$, in particular $A_{i,i-1}=k_i(p_in_i+q_im_i)$,
\item[] $P_i= \left\lbrace \lambda_1 X_{i-1}+ \lambda_2 X_i \ : \ 0 < \lambda_i < 1, \ i=1,2 \right\rbrace$, $1 \leq i \leq r$,
\item[] $L_i=] 0, 2X_i [$, $1 \leq i \leq r-1$
\item[] $\phi_i$ the linear map which take the value $1$ on $X_{i-1}$ and $X_i$, for $1 \leq i \leq r$, 
\item[] $\Delta_i= p_{i+1}q_{i}-p_{i}q_{i+1}$, $1 \leq i \leq r$.
\end{enumerate}

Let call $\Gamma_{ND}(p_1,q_1,k_1;\cdots;p_r,q_r,k_r)$ such a diagram.

\begin{remark}\rm
Due to the local situation we know that $\Delta_i>0$. We can permit $k_i=0$ by puting $\Gamma_{ND}(p_1,q_1,k_1;\cdots;p_r,q_r,k_r)=\Gamma_{ND}(p_1,q_1,k_1;\cdots;\hat{p}_i,\hat{q}_i,\hat{k}_i;\cdots;p_r,q_r,k_r)$ where the hat means that the term is omitted.
\end{remark}

From the previous sections we have the following well known equalities.

\begin{proposition}
The spectrum of $f$ is
$$
Sp(f)=\sum_{i=1}^r \sum_{(m,n)\in P_{i}\cap\N^2} (1-\phi_i(m,n)) + \sum_{i=1}^{r-1} \sum_{(m,n)\in L_i\cap\N^2} (1-\phi_i(m,n)).
$$
and the Milnor number of $f$ is
$$
\mu(f)=A_{0,-1}+A_{1,0}+\cdots+A_{r+1,r}+1.
$$
\end{proposition}

In \cite{Thomas Brelivet 1} it has been proved that we have the following theorem.

\begin{theorem}
We have
$$
6S-\mu\alpha_\mu= - \sum_{i=1}^{r-1}E_i \Delta_i,
$$
where 
$$
\alpha_\mu=1-\phi_{i_0}(1,1)
$$
is the maximal spectral value with $i_0$ such that $(1,1)$ is in the parallelogram generated by $X_{i_0-1}, X_{i_0}$ and
$$
E_i=
\left\lbrace
\begin{array}{ll}
\left( \sum_{k=-1}^{i_0-1} A_{k+1,k}(n_i-m_i)+d_i^2-m_i \right) \frac{k_ik_{i+1}}{A_{i,i-1}A_{i+1,i}}, \text{ if } 1 \leq i < i_0, \\
\left( \sum_{k=i_0}^{r} A_{k+1,k}(m_i-n_i)+d_i^2-n_i \right) \frac{k_ik_{i+1}}{A_{i,i-1}A_{i+1,i}}, \text{ if } i_0 \leq i \leq r.
\end{array}
\right.
$$
\end{theorem}

\begin{remark}\rm
The previous formula is true for every $i \in \{ 1, \ldots, r \}$, we use it for $i_0$ in order to prove the Hertling conjecture.
\end{remark}

The quantities $E_i$ are strictly positive so we deduce the following corollary.

\begin{corollary}
If $f$ is a germ Newton non-degenerated and commode then the conjecture of Hertling is true for $f$. Moreover, we have an equality if and only if $f$ is a positive deformation of a quasi-homogeneous polynomial defining an isolated singularity.
\end{corollary}

\begin{notation}\rm\label{remark : on the ci}
To simplify the computations we need to introduce:
$$
F_i=
\left\lbrace
\begin{array}{ll}
\sum_{k=-1}^{i_0-1} A_{k+1,k}(n_i-m_i)+d_i^2-m_i \text{ if } 1 \leq i < i_0, \\
\sum_{k=i_0}^{r} A_{k+1,k}(m_i-n_i)+d_i^2-n_i \text{ if } i_0 \leq i \leq r,
\end{array}
\right.
$$
and
$$
C_i=\frac{A_{i,i-1}+A_{i+1,i}-A_{i+1,i-1}}{A_{i,i-1}A_{i+1,i}}=-\frac{1}{(p_in_i+q_im_i)(p_{i+1}n_{i+1}+q_{i+1}m_{i+1})}\Delta_i.
$$
The formula becomes
$$
6S-\mu\alpha_\mu=\sum_{i=1}^{r-1}F_iC_i.
$$
\end{notation}

\begin{remark}\rm\label{remark : on the phii}
We have the following usefull eqality
$$
\frac{(n_i-m_i)\Delta_i}{(p_in_i+q_im_i)(p_{i+1}n_{i+1}+q_{i+1}m_{i+1})}=
\phi_i(1,1)-\phi_{i+1}(1,1).
$$

\end{remark}

\section{Newton polygonal representation and additivity of the spectral pairs}\label{section : Additivity of the spectral pairs}

Let $f : (\C^2,0) \rightarrow (\C,0)$ be an isolated singularity of curve.
From the Newton Puiseux algorithm, we can compute the \splice diagram of an isolated plane curve singularity. This algorithm follows a tree called here the polygon tree $\PP_f$ of $f$, $0$ being the root of $\PP_f$ and gives a splice diagram $\Gamma_f$. Each vertex of this tree correspond a polygon. 

In the following, without loss of generality we can supposse that we can use the representation of a \splice diagram given by the Figure \ref{Figure : general splice diagram} representing $f$.

\begin{figure}[H]
  \centering
  \input{GeneralSpliceDiagram.pstex_t}
  \caption{A general \splice diagram}\label{Figure : general splice diagram}
\end{figure}

\begin{definition}
Let $\mathcal{P}^a$ be the $\Z$-free module generated by the abstract Newton polygons $\Gamma_{ND}(p_1,q_1,k_1;\cdots;p_r,q_r,k_r)$ where $k_1, \cdots, k_r$ are positive integers and $(p_1,q_1),\cdots, (p_r,q_r)$ are positive coprime integers. 

Let $\sppa : \mathcal{P}^a \longrightarrow \Z^{(\Q)}$ the morphism such that $\sppa(\Gamma_{ND}(p_1,q_1,k_1;\cdots;p_r,q_r,k_r))=\spp(\Gamma_{ND}(p_1,q_1,k_1;\cdots;p_r,q_r,k_r))$.
\end{definition}

\begin{figure}[H]
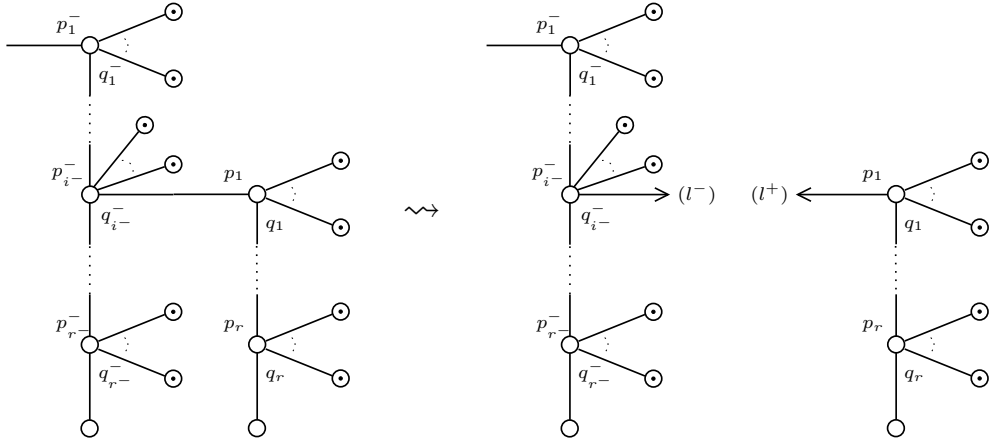

  \centering
  \input{consecutiveverticalparts.pstex_t}
  \quad
  \input{consecutiveverticalparts2.pstex_t}
  \caption{Two consecutive vertical parts of the \splice diagram}\label{Figure : Two consecutive vertical parts of the splice diagram}
\end{figure}

Cut the horizontal edges (those are not vertical) and put arrows with multiplicities such that the multiplicities of each vertex don't change. See Figures \ref{Figure : Two consecutive vertical parts of the splice diagram}. 

\begin{figure}[H]
  \centering
  \input{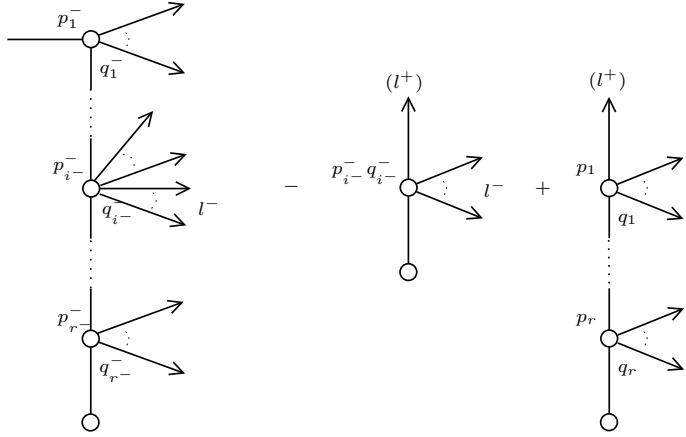}
  \caption{Cutted horizontal edges}\label{Figure: avantverticalnewtonpolygon}
\end{figure}

Then replace horizontal arrows with multiplicity $l$ by $l$ arrows of multiplicty $1$ ans substract by a new diagram in order to obtain polygonal diagrams. See Figure \ref{Figure: avantverticalnewtonpolygon}. 

\begin{figure}[H]
  \centering
  \input{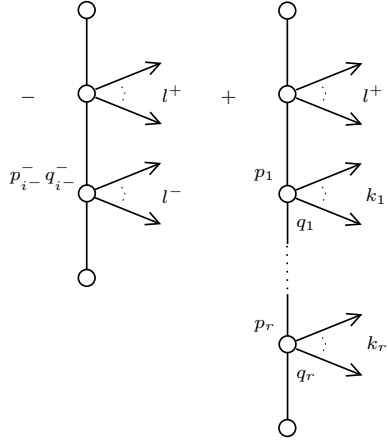}
  \caption{Vertical parts in $\mathcal{P}^a$}\label{Figure: verticalnewtonpolygon}
\end{figure}

If necessary ($l^+\not= 0$) complete the diagrams in order to have commode polygonal diagrams. See Figure \ref{Figure: verticalnewtonpolygon} for the last two parts.

In this way we have construct an element of $\P^a$: 
$$
\Psi(\Gamma_f)=\Psi_0(\Gamma_f)+\sum_{w \in \PP_f} \Psi_v(\Gamma_f)
$$

where $\Psi_0(\Gamma_f)=\Gamma_{ND}(p^0_1,q^0_1,k^0_1;\cdots;p^0_r,q^0_r,k^0_r)$ is the Newton polygon of the germ, 
$$
\Psi_w(\Gamma_f)=\Psi_w^+(\Gamma_f)-\Psi_w^-(\Gamma_f),
$$ 
$$
\Psi_w^+(\Gamma_f)=\Gamma(1,1,l^+;p_1,q_1,k_1;\cdots;p_r,q_r,k_r)
$$ 
and 
$$
\Psi_w^-(\Gamma_f)=\Gamma(1,1,l^+;p^-_{i^-}q^-_{i^-},1,l^-).
$$

\begin{theorem}
From $\Psi$, we get an equivalence relation $\sim$ on $\mathcal{P}^a$ and a map 
$$
\{ f \in \C \{ X,Y \} : f \text{ defines an isolated singularity}\} \rightarrow \mathcal{P}^a/\sim
$$
such that the composition with the morphism $\sppa$ gives the spectral pairs of $f$.
\end{theorem}

\begin{question}
What is $\ker(\sppa)$?
Do we have such a factorisation in higher dimensions?
\end{question}

\section{The maximal spectral value}\label{section : The maximal spectral value}

Recall that to $\Gamma(m,n,p,q,\l_1,\ldots,\l_a)$ we have associated
$$
\begin{array}{ccccc}
\phi & : & \R^2  & \rightarrow & \R \\
     &   & (x,y) & \mapsto     & \frac{qx+py}{qm+pn}.
\end{array}
$$
In the non-degenerated case the maximal spectral value is given by $1-\phi(1,1)$ where the $\phi$ corresponds to the parallelogram wich contains the point $(1,1)$ so it is natural to put the following definition.

\begin{definition}
We call $1-\phi(1,1)$ the virtual spectral value associated to $\Gamma(m,n,p,q,\l_1,\ldots,\l_a)$.
\end{definition}

We get a majoration of the spectral values (those of multiplicity non zero).

\begin{proposition}
We have
$$
1-\phi(1,1) \geq \alpha,
$$
for all spectral value $\alpha$ of $\Gamma_P(m,n,p,q,\l_1,\ldots,\l_a)$.
The virtual spectral value $1-\phi(1,1)$ is a spectral value if and only if $(1,1) \in \{ \lambda_0(m-p\l,n+q\l)+\lambda_1(m,n) : 0 \leq \lambda_0,\lambda_1 < 1 \}$. If it is not the case, then $\phi(1,1)$ is a spectral value.
\end{proposition}

We will show that this majoration is enough for our computations.

Now we can define the virtual spectral value associated to a vertex $v$ of the \splice diagram of the germ $f$.
Consider a vertex $v$ of the \splice diagram $\Gamma_f$ and the diagram associated to the vertex $v$ obtained by cuting all the edges around $v$. Then we get a diagram of the type $\Gamma(m,n,p,q,\l_1,\ldots,\l_a)$.

\begin{definition}
The virtual maximal spectral value of the vertex $v$ is the virtual spectral value associated to the previous diagram.
\end{definition}

To study the maximal spectral value of $\Gamma_f$ it is usefull to study the variation in the diagram of the virtual spectral values.

\begin{figure}[H]
  \centering
  \input{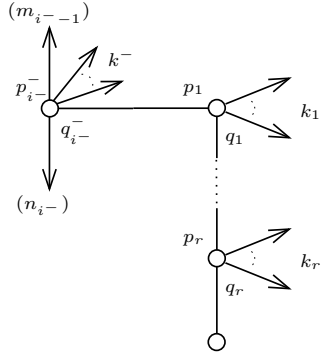}
  \caption{Attached vertical part}\label{Figure: Attached vertical part}
\end{figure}

Consider the part of the \splice diagram given by the Figure \ref{Figure: Attached vertical part}. 

Let 
$$
\l^+=p_{i^-}^-n_{i^-}+q_{i^-}^-m_{i^--1}+p_{i^-}^-q_{i^-}^-k^-
$$

and

$$
\l^-=q_1k_1+\cdots+q_rk_r.
$$

\begin{proposition}\label{Proposition : maximal vertical spectral value}
Consider an edge of the diagram of the Figure \ref{Figure: Attached vertical part}. Cut the edge and get two arrows with multiplicities. Then
\begin{enumerate}
\item[1--] if the edge is vertical the arrow with maximal multiplicity gives the direction where we have to go to get the virtual maximal spectral value,
\item[2--] if the edge is horizontal 
and $p_{i^-}^-, q_{i^-}^->1$ or $\l^+ > \l^-$ then the virtual maximal spectral value is given by the vertex of index $i^-$. 
\end{enumerate}
\end{proposition}

\proof
The first assertion is given for the vertical part (non-degenerate case) of the diagram by the remark \ref{remark : on the phii}.

The second assertion 
is given by the following computation.

Let 
$$
\alpha_{i^-}^-=1-\frac{p_{i^-}^-+q_{i^-}^-}{\l^++p_{i^-}^-q_{i^-}^-\l^-}
$$
and
$$
\alpha_j=1-\frac{p_j+q_j}{q_j(\l^++p_1k_1+\cdots+p_jk_j)+p_j(q_{j+1}k_{j+1}+\cdots+q_rk_r)}.
$$

We have to show that $\alpha_{i^-}^->\alpha_j$ for $j$ from $1$ to $r$.

We have (see section \ref{section : A formula for the variance of the spectrum} for the definition of $\Delta_{jb}$)

\bigskip

$q_j(\l^++p_1k_1+\cdots+p_jk_j)+p_j(q_{j+1}k_{j+1}+\cdots+q_rk_r) $
\begin{eqnarray*}
&= &q_j(\l^++p_1k_1+\cdots+p_jk_j)+p_j(\l^--q_1k_1-\cdots-q_jk_j) \\
&=&p_j\l^-+q_j\l^++\sum_{b=1}^{j-1}(p_bq_j-p_jq_b)k_b \\
&=&p_j\l^-+q_j\l^+-\sum_{b=1}^{j-1}\Delta_{jb}k_b
\end{eqnarray*}

the numerator of $\alpha_{i^-}^- - \alpha_j$ is
$$
[p_j+q_j-q_j(\pmm+\qmm)]\l^++[\pmm\qmm(p_j+q_j)-p_j(\pmm+\qmm)]\l^-+(\pmm+\qmm)\sum_{b=1}^{j-1} \Delta_{jb}k_b,
$$ 
and
$$
p_j+q_j-q_j(\pmm+\qmm)=q_j(\pmm-1)(\qmm-1)+(p_j-q_j\pmm\qmm),
$$
$$
\pmm\qmm(p_j+q_j)-p_j(\pmm+\qmm)=p_j(\pmm-1)(\qmm-1)-\Delta_0=p_j((\pmm-1)(\qmm-1)-1)+\pmm\qmm q_j.
$$
The numerator can be rewrite as
$$
(q_j\l^++p_j\l^-)(\pmm-1)(\qmm-1)+(\l^+-\l^-)\Delta_0+(\pmm+\qmm)\sum_{b=1}^{j-1}\Delta_{jb}k_b
$$
so it is positive if $\pmm,\qmm>1$ or $\l^+>\l^-$.

\endproof

{\bf Condition (Hl):} For each vertex of the diagram all the horizontal edges attached to a vertex except perhaps one are such that $\l^- \geq \l^+$.

\medskip

\begin{figure}[H]
\begin{center}
\input{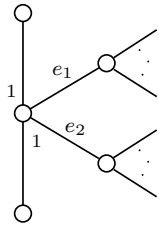}
\caption{Neighborhood of $v$}\label{Figure : twohl}
\end{center}
\end{figure}

Suppose that there is two horizontal edges $e_1, e_2$ attached (from the right) to the vertex $v \in V$ with decoration $(p,q)$ such that $\l_1^- \geq \l_1^+$ and $\l_2^- \geq \l_2^+$. By construction we have $\l_1^+=pq\l_2^-+a_1$, $\l_2^+=pq\l_1^-+a_2$, with $a_1,a_2 \geq 0$. Then $\l_1^+ \geq \l_2^-$ and $\l_2^+ \geq \l_1^-$. From this we deduce that $\l_1^+=\l_1^-=\l_2^+=\l_2^-$, $p=q=1$ and the diagram is in the form given by the Figure \ref{Figure : twohl}.

We can eliminate the two vertices of valence $1$ and the edge $e_1$ and $e_2$ become vertical. These operations don't change the spectrum and by this way we can suppose that we are in the situation of the condition (Hl).

\medskip

Now we can look for the maximal virtual spectral value of $\Gamma_f$.

\medskip

The point $(1,1)$ give us the vertex we must consider in the first Newton polygon (corresponding to the vertex $0 \in \PP_f$). Let $i_0^0$ be the index of this vertex.
Consider the horizontal edges attached to it. There is at most one edge such that $\l^- \geq \l^+$. If there is no such edge then the algorithm stops. If not consider the vertex attached to the right of the edge. It gives a new polygon (so a new vertex of $\PP_f$) and we consider the new vertex corresponding given by the point $(1,1)$. If $p_{i_0^0}$ and $q_{i_0^0}$ are different to $1$ then we continue the algorithm with the new vertex instead of the vertex of index $i_0$. If not the we exchange the horizontal edge with the vertical one of decoration $1$ and we continue the algorithm.

If we do again the algorithm then we get a sequence of vertices (and the tree is not modified!) $v_{i_0^0}, \ldots, v_{i_0^n}$. Let $\alpha_{i_0^0}, \ldots, \alpha_{i_0^n}$ the corresponding virtual spectral values and $(p_{i_0^0},q_{i_0^0}), \ldots,$ $(p_{i_0^n},q_{i_0^n})$ the decorations attached to these vertices.
The algorithm tells us that the decorations $p_{i_0^0},q_{i_0^0}, \ldots, p_{i_0^n},q_{i_0^n}$ are all different to $1$ and $\alpha_{i_0^0} > \cdots > \alpha_{i_0^n}$ (case 2 of the proposition).
Take now any vertex $v$ of the tree $\Gamma_f$ with virtual spectral value $\alpha_v$ and consider the geodesic from $v$ to $v_{i_0^0}$. Until we have not reach a vertex of the sequence $v_{i_0^0}, \ldots, v_{i_0^n}$, in order to find the maximal virtual spectral value we follow the direction given by the arrow of maximal multiplicity. From the Proposition \ref{Proposition : maximal vertical spectral value} the virtual spectral values are growing. Call $v_{i_0}^k$ the first vertex of the sequence we met. To conclude, we have $\alpha_v \leq \alpha_{i_0^k} \leq \alpha_{i_0^0}$.

\begin{proposition}
The virtual spectral value 
$$
\alpha_{i_0^0}
$$
is the maximal virtual spectral value.
\end{proposition}

\begin{proposition}
The maximal virtual spectral value is the maximal spectral value.
\end{proposition}

\proof
We just have to check that the virtual spectral value $\alpha_{i_0^0}$ is a spectral value.

From the choice of $i_0^0$ the value $\alpha_{i_0^0}$ has multiplicity $1$ in 
$$
\sp(\Gamma(m_{i_0^0+1},n_{i_0^0},p_{i_0^0},q_{i_0^0},k_{i_0^0})).
$$ 
The spectral values of $\sp(\Gamma(\l^+_{i_0^0}+p_{i_0^0}q_{i_0^0}\l_{i_0^0}^-,0,p_{i_0^0}q_{i_0^0},1,\l_{i_0^0}^-))$ are of the form 
$$
1-\frac{x+p_{i_0^0}q_{i_0^0}y}{\l_{i_0^0}^++p_{i_0^0}q_{i_0^0}\l_{i_0^0}^-}, \quad x,y \geq 1
$$
and from the fact that $p_{i_0^0}>1, q_{i_0^0}>1$ it is strictly greater than 
$$
\alpha_{i_0^0}=1-\frac{p_{i_0^0}+q_{i_0^0}}{p_{i_0^0}n_{i_0^0}+q_{i_0^0}m_{i_0^0}}
=1-\frac{p_{i_0^0}+q_{i_0^0}}{\l_{i_0^0}^++p_{i_0^0}q_{i_0^0}\l_{i_0^0}^-}.
$$
\endproof

\begin{theorem}\label{Theorem : The spectral maximal value}
Let $f : (\C^2,0) \rightarrow (\C,0)$ be an analytic germ defining an isolated singularity of curve. Then there exists a germ $g : (\C^2,0) \rightarrow (\C,0)$ topologically equivalent to $f$ such that the maximal spectral value is given by the point $(1,1)$ in the plane of the Newton polygon of $g$. Furthermore the multiplicity of the maximal (or minimal) spectral value is one.
\end{theorem}

\begin{remark}\rm
In general the germ $g$ is not Newton non-degenerated.
\end{remark}

\section{A formula for the variance of the spectrum, Hertling conjecture}\label{section : A formula for the variance of the spectrum}

From the previous section we can suppose that the maximal spectral value of $f$ is given by the Newton polygon.

We have seen that we have a decomposition of the splice diagram as an abstract sum of polygons,
$\Psi(\Gamma_f)=\Psi_0(\Gamma_f)+\sum_{w \in \PP_f}(\Psi_w^+(\Gamma_f)-\Psi_w^-(\Gamma_f))$.

Let $w \in \PP_f$ different of the root $0$ and
\begin{enumerate}
\item[] $\mu_0$ the Milnor number of $\Psi_0(\Gamma_f)$,
\item[] $\mu^+_w$ the Milnor number of $\Psi_w^+(\Gamma_f)$,
\item[] $\mu^-_w$ the Milnor number of $\Psi_w^-(\Gamma_f)$,
\item[] $\mu_w=\mu^+_w-\mu^-_w$,
\item[] $S^+_w$ the sum of the squares of the spectral values of $\Psi_w^+(\Gamma_f)$,
\item[] $S^-_w$ the sum of the squares of the spectral values of $\Psi_w^-(\Gamma_f)$,
\item[] $S_w=S^+_w-S^-_w$,
\item[] $\alpha_0$ the maximal spectral value of $\Psi_0(\Gamma_f)$ (and of $\Gamma_f$),
\item[] $\l^-=q_1k_1+\cdots + q_rk_r$, and $\l^+$ is the multiplicity of the left arrow, 
\item[] $\alpha^+_w=1-\frac{p_1+q_1}{p_1\l^-+q_1\l^+}$, 
\item[] $\alpha^-_w=1-\frac{2}{\l^++\l^-}$, 
\item[] $p_0=p_{i^-}^-q_{i^-}^-$, $q_0=1$, $\alpha_0^-=1-\frac{p_0+q_0}{p_0\l^-+q_0\l^+}$,
\item[] $\Delta_{i,j}=p_iq_j-p_jq_i$, for $0 \leq i, j \leq r$,
\item[] $\Delta_i=p_{i+1}q_{i}-p_{i}q_{i+1}$, $0 \leq i \leq r-1$.
\end{enumerate}

From the additivity of the Milnor number and the sum of the squares of the spectral values we get the following lemma. 

\begin{lemma}
$$
6S-\mu\alpha_\mu=
6S_0-\mu_0\alpha_0 +
\sum_{w \in \PP_f} \left[ (6S^+_w-\mu^+_w\alpha^-_w)-(6S^-_w-\mu^-_w\alpha^-_w)+ (\mu^+_w-\mu^-_w)(\alpha^-_w-\alpha_0) \right].
$$
\end{lemma}

We already know from the Newton and non-degenerated case that $6S_0-\mu_0\alpha_0$ is a linear combination with negative coefficients of the determinants of the \splice diagram.
We will show that if we fix a vertex $w$ of $\PP_f$, then $(6S^+_w-\mu^+_w\alpha^-_w)-(6S^-_w-\mu^-_w\alpha^-_w)$ and $\mu^+_w-\mu^-_w$ are linear combination of the determinant correponding to the vertex $w$. 
Finaly we will show that $6S-\mu\alpha_\mu$ is a linear combination of the determinants with negative coefficients. From this result we deduce immediately the Hertling conjecture.

From the non-degenerated and commode case we have
\begin{enumerate}
\item[] $6S^+-\mu^+\alpha^-= F_0^+C_0^++\cdots + F_{r-1}^+ C_{r-1}^+$,
\item[] $6S^--\mu^-\alpha^-=F_0^-C_0^-$,
\end{enumerate}
where we have forgotten $w$ to simplify notations.

Suppose that $l^+ \geq l^-$ then we know that the signs of the $F_i$ are positive. In the case of $l^- > l^+$, the same type of computation will work.

\begin{lemma}
\begin{multline*}
F_0^+C_0^+ - F_0^-C_0^-=
\left[ d^2-\l^- + (\l^+ + \pp\qq \l^-)(\l^+-\l^-)(\l^--1) \right] (C_0^+-C_0^-)\\
+(\mu^+-\mu^-)C_0^+ (\l^+-\l^-).
\end{multline*}
\end{lemma}

\proof
From the section \ref{section: the non-degenerated and commode case} we have:
$$
\begin{array}{rcl}
F^+_0 & = & (A^+_{2,1}+A^+_{3,2}+\cdots+A^+_{r+2,r+1})(\l^+-\l^-)+d^2-\l^- \\
      & = & (\mu^+ - (\l^+-1)(\l^++\l^-)-1)(\l^+-\l^-) +d^2-\l^-
\end{array}
$$

and

$$
\begin{array}{rcl}
F^-_0 & = & (A^-_{2,1}+A^-_{3,2})(m_1^--n_1^-)+d^2-\l^- \\
      & = & (\mu^- - (\l^+-1)(\l^++\l^-)-1)(\l^+-\l^-) +d^2-\l^-.
\end{array}
$$
So
\begin{multline*}
F_0^+C_0^+ - F_0^-C_0^-=
\left[ d^2-\l^--((\l^+-1)(\l^-+\l^+)+1)(\l^+-\l^-)+\mu^-(\l^+-\l^-) \right] (C_0^+-C_0^-)\\
+(\mu^+-\mu^-)C_0^+ (\l^+-\l^-).
\end{multline*}

\endproof

\begin{lemma}
$$
C_0^+-C_0^-= \frac{-\Delta_0}{(q_1\l^++p_1\l^-)(\l^++\pp\qq \l^-)}.
$$
\end{lemma}

\proof
From the definitions we have:
$$
C_0^+=\frac{q_1-p_1}{(\l^++\l^-)(q_1\l^++p_1\l^-)}
$$
and
$$
C_0^- = \frac{1-\pp\qq}{(\l^++\l^-)(\l^++\pp\qq \l^-)}.
$$
\endproof

\begin{lemma}
$$
\mu^+-\mu^-=\sum_{t=0}^{r-1} \frac{n_{t}(n_{t}-1)}{q_tq_{t+1}} \Delta_t.
$$
where $n_0=\l^-$.
\end{lemma}

\proof
From the definitions we have:
\begin{multline*}
\mu^+=(\l^++\l^-)(\l^+-1)
+ \sum_{i=1}^{r} \left[ p_i(q_ik_i + \cdots + q_rk_r) + q_i(\l^++p_1k_1+ \cdots + p_{i-1}k_{i-1}) \right] k_i \\
-(\l^++p_1k_1+\cdots+p_rk_r)+1.
\end{multline*}
and
$$
\mu^-=(\l^++\l^-)(\l^+-1)
+ (\pp\qq \l^-+\l^+)(\l^--1)+1.
$$
So
\begin{multline*}
\mu^+-\mu^-=-\pp\qq \l^-(\l^--1)+\sum_{i=1}^r
\left[ p_i(q_ik_i + \cdots + q_rk_r) + q_i(p_1k_1+ \cdots + p_{i-1}k_{i-1}) \right] k_i  \\
-(p_1k_1+\cdots+p_rk_r),
\end{multline*}
\begin{multline*}
\mu^+-\mu^-=-\pp\qq \l^-(\l^--1)
+ \left( \sum_{i=1}^r p_ik_i \right) (\l^--1) \\
+ \sum_{i=1}^r \left[ (p_1q_i-p_iq_1) k_1 + \cdots + (p_{i-1}q_i - p_iq_{i-1})k_{i-1} \right] k_i.
\end{multline*}
From the equality
$$
\sum_{i=1}^rp_ik_i=\frac{p_1}{q_1}\l^- + \sum_{i=2}^r \frac{p_iq_1-p_1q_i}{q_1}k_i.
$$
we get

$$
\begin{array}{rcl}
\mu^+-\mu^- & = &
\frac{p_1-\pp\qq q_1}{q_1} \l^-(\l^--1) + \sum_{i=2}^{r} \left( \frac{\l^--1}{q_1}(p_iq_1-p_1q_i)+ \sum_{j=1}^{i-1}(p_jq_i-p_iq_j)k_j\right)k_i. \\
 & = & \frac{\Delta_0}{q_1}\l^-(\l^--1) + \sum_{i=2}^r \left( \frac{\l^--1}{q_1} \Delta_{i,1} - \sum_{j=1}^{i-1} \Delta_{i,j} k_j \right)k_i.
\end{array}
$$
We also have
$$
\Delta_{i,j}= \sum_{t=j}^{i-1} \frac{q_iq_j}{q_tq_{t+1}} \Delta_t
$$
so
$$
\sum_{j=1}^{i-1} \sum_{t=j}^{i-1} \frac{q_iq_j}{q_tq_{t+1}}\Delta_tk_j = q_i\sum_{t=1}^{i-1} \frac{\Delta_t}{q_tq_{t+1}} \sum_{j=1}^{t}q_jk_j.
$$
We get the result from the equality:
$$
\begin{array}{rcl}
\mu^+-\mu^- & = & \frac{\Delta_0}{q_1}\l^-(\l^--1) + \sum_{i=2}^r \sum_{t=1}^{i-1} q_i \frac{\Delta_t}{q_tq_{t+1}} \left( n_{t}-1 \right) k_i, \\
            & = & \frac{\Delta_0}{q_1}\l^-(\l^--1) + \sum_{t=1}^{r-1} \sum_{i=t+1}^r q_i \frac{\Delta_t}{q_tq_{t+1}} \left( n_{t}-1 \right) k_i.
\end{array}
$$
\endproof

\begin{lemma}
$$
\alpha^+-\alpha^-=C_0^+(\l^+-\l^-).
$$
\end{lemma}

\proof
From the definitions.
\endproof

From the previous lemmas, we deduce:

\begin{proposition}
\begin{multline*}
(6S-\mu\alpha_\mu)_w =
- \left( \frac{d^2-\l^-}{(q_1\l^++p_1\l^-)(\l^++p_i^-q_i^-\l^-)} - \frac{\l^--1}{q_1}(\l^--1-\l^-\alpha_0) \right) \Delta_0 \\
- \sum_{t=1}^{r-1} \left( E_t^+ + \frac{n_{t}(n_{t}-1)}{q_tq_{t+1}}(\alpha_0-\alpha^+) \right) \Delta_t.
\end{multline*}
\end{proposition}

\proof
From the previous Lemmas we have:
\begin{multline*}
(6S-\mu\alpha_\mu)_w=-E_2^+\Delta_2-\cdots-E_{r-1}^+\Delta_{r-1} \\
+ \left[ d^2-\l^- + (\l^+ + \pp\qq \l^-)(\l^+-\l^-)(\l^--1) \right] (C_0^+-C_0^-) \\
+(\mu^+-\mu^-)(\alpha^+-\alpha_0).
\end{multline*}
So we have immediately the coefficient of $\Delta_t$ for $t >0$ and for $\Delta_0$ we have to remark that 
$$
\alpha^+ \frac{\l^-}{q_1}-\frac{\l^+-\l^-}{q_1\l^++p_1\l^-}=\frac{\l^--1}{q_1}.
$$

\endproof

Now it is easy to prouve that the Hertling conjecture is true in dimension $2$.

\begin{theorem}
There exists positive rational numbers $(E_e)_{e \in \ed}$ such that
$$
6S-\mu\alpha_\mu=-\sum_{e \in \ed}E_e \Delta_e.
$$
\end{theorem}

\proof
We just have to show that the first coefficient is positive.

We know that $\alpha_0$ is greater than 
$$
\alpha_i=1-\frac{\pp+\qq}{\l^++\pp\qq \l^-}
$$
so we only have to check the inequality for $\alpha_i$ instead of $\alpha_0$.

We have
\begin{multline*}
q_1(d^2-\l^-)+(q_1\l^++p_1\l^-)(\l^++\pp\qq \l^--\l^-(p_i^-+q_i^-))(\l^--1) = \\
(d^2-1)q_1 + \left[ (q_1\l^++p_1\l^-)(\l^++((p_i^--1)(q_i^--1)-1)\l^-)-q_1 \right] (\l^--1)
\end{multline*}
If $\pp, \qq>1$ then the last expression is positive. If $\pp$ or $\qq$ is equal to $1$ then $\l^+>\l^-$ then the last expression is also positive.
\endproof

\begin{corollary}\label{corollary : Hertling conjecture}
If $f : (\C^2,0) \rightarrow (\C,0)$ is an isolated singularity then the conjecture of Hertling is true for $f$. Moreover, we have an equality if and only if $f$ is a positive deformation of a quasi-homogeneous polynomial defining an isolated singularity.
\end{corollary}

\begin{flushleft}
Thomas Brélivet \\
Departamento de Álgebra, Geometría y Topología \\
Facultad de ciencias \\
Universidad de Valladolid \\
47005 Valladolid \\
España \\
email: brelivet@agt.uva.es \\
\end{flushleft}


\begin{thebibliography}{9}

\bibitem[AGV]{AGV}  V.I. Arnold, S.M. Gusein-Zade and A.N. Varchenko, 
  {\it Singularities of differentiable maps}, 
       vol 1,2, Birkh\"auser, Boston, 1988.

\bibitem[B1]{Thomas Brelivet 1} T. Br\'elivet, {\it Variance of the spectral number and Newton polygons}, Bull. Sci. Math. 126 (2002), no. 4, 332--342.

\bibitem[B2]{Thomas Brelivet These} T. Br\'elivet, {\it Topologie des polyn\^omes, spectre et variance du spectre}, Th\`ese Universit\'e Bordeaux 1, 2002.

\bibitem[B3]{Thomas Brelivet 3} T. Br\'elivet, {\it Sur les pairs spectrales de polyn\^omes \`a deux variables}, Preprint.

\bibitem[C]{C} P. Cassou-Nogu\`es, {\it Entrelacs toriques it\'er\'es et int\'egrales associ\'ees \`a une courbe plane}, S\'eminaire de Th\'eorie des Nombres, Bordeaux  2 (1990), 273-331. 

\bibitem[CP]{CP} P. Cassou-Nogu\`es, A. Ploski, {\it Introduction to Algebraic Plane Curve Singularities}, in preparation. 

\bibitem[Di]{Di} A. Dimca, {\it Monodromy and Hodge theory of regular functions}, Proceedings of the Summer Institute on Singularities, Newton Institute, Cambridge 2000.

\bibitem[EN]{Eisenbud Neumann} D. Eisenbud, W. Neumann, {\it Three-dimensional link theory and invariants of plane curve singularities}, Annals of Mathematics Studies, 110, Princeton University Press, 1985.

\bibitem[H]{H} C. Hertling, {\it Frobenius manifolds and variance of the spectral numbers}, Proceedings of the Summer Institute on Singularities, Newton Institute, Cambridge 2000, see also Preprint math.CV/0007187.

\bibitem[KV]{Khovanskii Varchenko} A. G. Khovanski\u\i, A. N. Varchenko, {\it Asymptotic behavior of integrals over vanishing cycles and the Newton polyhedron}, Dokl. Akad. Nauk SSSR 283 (1985), no.~3, 521--525.

\bibitem[N]{Neumann} W.D. Neumann, {\it Topology of hypersurface singularities}, Preprint. \url{http://www.math.columbia.edu/~neumann/preprints/kaehler1.ps}

\bibitem[S1]{S1} M. Saito, {\it Period mapping via Brieskorn modules}, Bull. Soc. math. France 119 (1991), 141--171.

\bibitem[S2]{S2} M. Saito, {\it Exponents of an irreducible plane curve singularity}, Preprint math.AG/ 0009133. 

\bibitem[SS]{SS} J. Scherk and J. Steenbrink, {\it On the mixed Hodge structure on the cohomology of the Milnor fiber}, Math. Ann. 271 (1985), 641-665.

\bibitem[SSS]{SSS} R. Schrauwen, J. Steenbrink, J. Stevens, {\it Spectral pairs and the topology of curves singularities}, Proceedings of Symposia in Pures Mathematics, Volume 53 (1991), 305-328. 

\bibitem[St1]{St1} J. Steenbrink, {\it Mixed Hodge structures on the vanishing cohomology}, in P. Holm ed.: Real and Complex Singularities, Oslo 1976.

\bibitem[St2]{St2} J. Steenbrink, {\it Applications of Hodge theory to singularities}, Proc. I.C.M. Kyoto (1990), 559-576.

\bibitem[V]{V}  A. N. Varchenko, {\it The asymptotics of holomorphic forms determine a mixed Hodge structure}, Sov. Math. Dokl. 22(1980), 248-252).

\end{thebibliography}
\end{document}